\documentclass[12pt]{article}
\usepackage{graphicx}
\begin{document}
\title{Graph Entropy, Network Coding and Guessing games}
\author{ S\o ren Riis} 
\date{November 25, 2007}
\maketitle
\begin{abstract}
We introduce  {\it the (private) entropy of a directed graph} (in a new network coding sense) as well as a number of related concepts. We show that the entropy of a directed graph is identical to its guessing number and can be bounded from below with  the number of vertices minus the size of the graph's shortest index code.  We show that the Network Coding solvability of each specific multiple unicast network is completely determined by the entropy (as well as by the shortest index code) of the directed graph that occur by identifying each source node with each corresponding target node.

Shannon's information inequalities can be used to calculate upper bounds on a graph's entropy
as well as calculating the size of the minimal index code.  Recently, a number of new families of 
so-called non-shannon-type information inequalities have been discovered. It 
has been shown that there exist communication networks with a capacity strictly 
ess than required for solvability, but where this fact cannot be derived using 
Shannon's classical information inequalities. Based on this result we show that 
there exist graphs with an entropy  that cannot be 
calculated using only Shannon's classical information inequalities, and show 
that better estimate can be obtained  by use of certain non-shannon-type 
information inequalities.

%We notice that a long standing open conjecture by Valiant can be restated in terms of the behavior the minimal index code of shifts of labeled sparse graphs.  We also  notice that another long standing open question - the matrix-transposition problem (as well a related problems concerning improving the efficiency of  the discrete fourier transformation) - can be viewed as a reverse engineering problem for multiple unicast networks in Network Coding. 
\end{abstract}

\noindent
{\bf Category} E.4 {\it Graph Theory, Information Theory, Network Coding and Circuit Complexity}
%A category including the fourth, optional field follows...
%\category{D.2.8}{Software Engineering}{Metrics}[complexity measures, performance measures]
%\terms{Theory}
%\newpage
%\keywords{} % NOT required for Proceedings

\section{Introduction}
\subsection{Main results}
Informally, one of the key problems in proving lower bounds on circuits is to identify and prove that certain Circuit topologies (e.g. small circuits) provide information bottlenecks for information flows. 
We introduce the notion of (private/public) graph entropy and show that this notion in a very precise way captures such information bottlenecks in the context of communication networks.  More specifically we show that a given communication network $N$ with $k$ source and $k$ corresponding target nodes provide an information bottleneck (i.e. is unsolvable) if and only if the entropy (or public entropy) of $G_N$ is strictly less than $k$. In the seminal works \cite{ZY97,ZY97A} and  \cite{DFZ07} it was shown that Shannon's information inequalities (commonly known as Shannon's laws of information theory)
are in general insufficient to identify specific information bottlenecks in communication flow problems. 
Relying on this result we will show that a similar type of result is valid in the context of graphs entropy.
By combining Theorem 1, Theorem 9 as well as Theorem 10, and modifying  the acyclic multiple unicast version of the Vamos graph introduced in \cite{DFZ07} (by passing to the corresponding line graph and by identifying each source node with each corresponding target node) we show that the resulting graph  has an entropy that cannot be calculated using Shannon's Classical information inequalities and that better bounds can be achieved by use of Zhang and Young's non-shannon-information inequalities. \footnote{later and joint with Sun Yun, we found - partly by extensive computer searches - somewhat smaller and simpler examples }

In the paper we introduce a number of new (directed) graph parameters\footnote{throughout the paper graphs are always assumed to be directed}. These include the (private/public) entropy,  the Shannon (private/public) entropy,  the Zhang-Young (public/ private) entropy  as well as the Doghtery-Freiling-Zeger (private/public) entropy.  The Concepts are in general different from (but linked to) graph parameters that have already been extensively studied in the literature.  
In general we drop the prefix "private" (when no confusion is possible) and refer to the (private) graph parameters as Entropy, S-entropy, ZY-entropy and DFZ-entropy.   

\subsection{Experimental results for small graphs}

Using computer calculations we  tested millions of (small) graphs and found that quite different notions
(different variants of guessing numbers and Graph Entropy) led to identical numbers on the vast majority of graphs\footnote{a more detailed account of these and many other experimental findings is being prepared in collaboration with Sun Yun}.  Key concepts (most of which we will introduce in this paper), like the Graph Entropy, the guessing number, the S-entropy and the ZY-entropy - gave identical results for the vast majority of graphs we tested. One interesting aspect is that although the different calculations usually lead to the same numbers for the graphs we tested by computer (using weeks of computer time), there appears to be a very sparse set of (small) graphs where the calculations leads to slightly different results. 
It should, however, be noticed that small graphs are atypical and it would be a mistake to expect that the computational findings for small graphs with $\leq 10$ nodes in general remains valid for large graphs.  To illustrate the problem, we noticed for example that most small oriented graphs  have entropy identical to the number of vertex minus the acyclic independence number (i.e. the number of vertices in the maximal induced acyclic subgraph). This  is,  we believe, not typical for example, for large complex oriented graphs (e.g. random tournaments).

One curious fact (we discovered experimentally), is that the entropy for the vast majority of (small) graphs is an integer. More specifically, we investigated (using computer calculations) many classes of graphs. The vast majority of "small" graphs (less than $10$ nodes) seem to have integer or half integer Entropy (identical to their S-entropy and ZY-entropy) \footnote{we found a few graphs with Shannon Entropy and ZY-entropy being a third integer}.

This is striking as there does not seem to be any obvious reason why graphs in general should have integer Entropy. If we chose for the complete graph $K_n$ for each edge an edge direction (getting a tournament on $n$ vertex) our computer calculations showed that there is $1$ graph of type $K_2$,  $2$ non-isomorphic graphs of type $K_3$, $4$ non-isomorphic graphs of type $K_4$, $12$ non-isomorphic graphs of type $K_5$, 
$56$ non-isomorphic graphs of type $K_6$, $456$ non-isomorphic graphs of type $K_7$ and $6880$
of type $K_8$ and $191536$ non-isomorphic graphs of type $K_9$ \footnote{we identified this sequence independently by use of the dictionary of integer sequences}.

Our computer calculations showed that all, but one, of the $531$ non-isomorphic graphs of type $K_j$ for $j \leq 7$ have integer Entropy. 
We did not investigate all graphs of type $K_8$ and $K_9$ \footnote{each graph required a couple of hours computer analysis}, but random checks suggest that the majority of these graphs have integer Entropy and that most remaining graphs has half integer Entropy. We do not know if this pattern changes for graphs of type $K_n$ for $n$ large.  What we do know - we notice this in section 7 - is that the entropy is an integer whenever the (private) entropy and the public guessing number (=number of 
vertices minus the smallest index code) is identical.

Another striking point is that the (Shannon) Entropy for the vast majority of (small) graphs we tested is identical to the (Shannon) Entropy of the dual graph where all edge directions have been reversed.
To appreciate this we noticed that the actual calculations needed to workout these identical Entropies, are often very different, and do not seem to mirror each other in any obvious sense.  

\subsection{Historical background}
In \cite{Val} based on \cite{Valorg} Valiant introduced the concept of a bi-parte graph representing $n$ input output pairs, being realized by $m$ common bits. It turns out that Valiant's concept of being realized by $m$ common bits, is mathematically equivalent to the graph $G$ naturally associated to the bi-parte graph having an index code of length $m$.  The notion of an index code of a graph (as well as a minimal index code) was first introduced in \cite{Tel98} and further result on this graph parameter was presented in \cite{Tel06} and \cite{Tel07}  \footnote{we recall the definition in section 7}.
To confuse matters further, the very same concept was  independently referred to in \cite{Riis05} and \cite{Riis07} as $G$ having information defect $m$. To avoid further confusion we will follow
\cite {Tel98,Tel06} and \cite{Tel07} and refer to the notion as the minimal index code. 

Given a graph $G=(V,E)$ with vertex set $\{0,1,2,\ldots,n-1\}$ we define for each $t \in \{0,1,2,\ldots,n-1\}$ the "shifted" graph $G^t=(V,E^t)$ with vertex set $V$ and with and edge $(v_1,v_2) \in E$
exactly when $(v_1, v_2 \oplus t) \in E^t$ where $v_2 \oplus t$ is calculated modulo $n$. This construction was essentially presented in \cite{Val18}.

Valiant's shift problem asks whether for each $\epsilon >0$ for all sufficiently large values of $n$,
each sparse graph $G=(V,E)$ with vertex set $\{0,1,2,\ldots,n-1\}$ and at most 
$n^{1+\epsilon}$ edges, can be shifted such that $G^t$ has minimal index code of length at least $m$.
If it could be shown to be valid for $m=O(\frac{n}{\log(\log(n))})$, then a significant breakthrough in Boolean circuit complexity would follow \cite{Pudlak97} \cite{Riis07}.  In \cite{Valorg} Valiant conjectured that $m <n/2$ was not achievable, and before \cite{Riis07} it was not even known if $m < n -o(n)$ could be achieved. However in \cite{Riis07} it was shown that in certain cases all cyclic shifts lead to graphs with a (surprisingly) short minimal index code of length less than $m=(n-n^\epsilon)/2$. Despite this result, Valiant's original common information approach remains viable.  The guessing number of a graph was introduced in \cite{Riis05} and  was later related to Valiant's shift problem \cite{Riis07}.

The notion of Guessing numbers \cite{Riis05} and the notion of (private) graph entropy has, in general, been guided by analysis of specific communication flow problems. Our notion of guessing number of a graph 
grew out of work on the matrix transposition problem that is another long-standing open problem in circuit complexity \cite{Vitnew}.  To illustrate some of the difficulties in solving the matrix transposition problem, in \cite{Riis04} it was shown that there exist network flow problems that can be only solved by use of non-linear boolean functions. More specifically
we constructed  a specific acyclic graph with $25$ input nodes and $25$ output nodes, such that any linear solution (over an alphabet of two letters) can only transmit $24$ of the messages correctly, while there exists a non-linear solution that correctly transmits all $25$ messages. In the same paper 
\cite{Riis04} we presented a general construction, that allowed to lift this (and similar results) that are valid for alphabets with two letters, to alphabets with any number of letters. However, for each size of alphabet the construction leads to a new (and typically larger) graph. In \cite{DFZ05} the authors showed that such an increase in graph size can be avoided and that there exists an  acyclic communication network that has a non-linear solution  over an alphabet of $4$-letters, but fails to have linear solution (even if the alphabet is organized as a vector space) over alphabets of any finite cardinality. The paper \cite{DFZ06} was also significant as the authors used arguments from information theory to show unsolvability of specific communication networks. 

In \cite{Leh} and \cite{med} the authors showed that a network might be solvable over all alphabets, but not scalar-linearly solvable over any alphabet.

In \cite{Tel07} the authors used a very interesting general construction (based on an earlier construction that Alon, used to answer a question about Shannon Capacity \cite{Alon}). They showed that there exist  graphs $G_{p,q}$ with no linear index code of less that $n^{1-o(1)}$ bits, but with a non-linear index code of $n^{o(1)}$ bits. It is still an open question if there exists solvable communication networks, that requires non-linear coding functions, and where the performance (e.g. the network capacity) gets substantially reduced if all coding functions are linear. 

In \cite{Riis06} a we constructed a specific acyclic graph that we showed constitutes what we will call a {\it one-way information bottleneck}. More specifically in \cite{Riis06} a concrete acyclic graph with $38$ source nodes and $38$ target nodes were constructed, such that the messages fed into the source nodes in general can be sent to their corresponding target nodes, while the "dual" network where all edge directions have been reversed and the role of source nodes and target nodes have been reversed - is unsolvable (only $37$ messages can be sent in that direction). Using the lifting method from 
\cite{Riis04} for any size $s \leq 2$ of alphabet, it is possible to construct a one-way communication network that is solvable over alphabets of size $s$ in one direction, but is unsolvable in the reverse direction.

The one-way information bottleneck constructed in \cite{Riis06} was improved in \cite{DFZ06A} where it was shown  that there is a network that has a solution in one direction, but is unsolvable in the opposite direction. In the unsolvable direction the network remains unsolvable even if the size of alphabet is allowed to have any finite size \footnote{in the solvable direction only certain sizes of the alphabet is allowed}. It should be pointed out that this network in general (also for large alphabets) is unsolvable in both directions. It is an open question if there exists a network that is a genuine one-way bottleneck, that is solvable for any sufficiently large alphabet in one direction, but is unsolvable for any choice of alphabet in the opposite direction. It is also an open question if more dramatic  one-way bottlenecks can be constructed (e.g. with a polynomial difference in performance). The graph in 
\cite{Tel07} is self-dual (in fact all edges are bidirectional), so this construction would need to be modified to work.
In \cite{rev} and independently in \cite{Riis06} is was shown that networks that are solvable by linear coding functions  never contain a one-way information bottleneck. 

\section{The Guessing Number of a graph}

The Guessing number of a (directed) graph was first introduced in \cite{Riis05} and \cite{Riis07}. To explain this concept consider the following situation (based on \cite{Riis04}) where $100$ players each has a die with $s$-sides (you may assume the die is a traditional die, i.e. that $s=6$). Assume each person rolls their die. No-one is allowed to look at their own die, however each person has access to the dice values of the remaining $99$ people. One way of presenting this is that each player has the value of their own die
stuck to their forehead. But, unlike the famous puzzle where the players can get additional information from other players' hesitation, in this game, all players have to "guess"- simultaneously and without any communication - what die value they have stuck to their forehead. 
{\it What is the probability that all  $100$ players correctly "guess" the value of their own die?} 
It turns out that this question is ill-posed since the probability depends on the general protocol adopted by the players. If, for example, each player assumes that the sum of all dice-values is divisible by $s$ and each player from this assumption deduces the value of their own die - all players are correct
when exactly the sum of the dice value is indeed divisible by $s$. Or equivalently if just one player is right (this happens with probability $\frac{1}{s}$) all players are right. 
Had the players guessed in a random uncoordinated fashion they would only all have been right
with probability $\frac{1}{s}^{100}$. Thus adopting the "$0$-mod $s$" guessing strategy 
the players are $s^{99}$ times more likely to be right, than if they guess randomly in an uncoordinated fashion.

This game can be modified and played on a directed graph $G=(V,E)$. The idea is that each node represents a player, and again each player throws a die with $s$-sides.  The value of each die is passed on along each edge. Thus the player at node $j$ has access to only the dice values corresponding to players in nodes $i$ where $(i,j)$ is an edge in the graph. The task is - as before - to find a guessing strategy that maximize the probability all $n=|V|$ players guess correctly the value of their die.

The example with the $100$ players corresponds to the complete graph $K_{100}$ (with bi-directed edges) on $100$ nodes. On $K_{100}$ the players can do $s^{99}$ times better than pure random guessing, which is why we say that the complete graph $K_{100}$ has guessing number $99$. In general:

\begin{description}

\item [Definition]: 
\it 
\noindent
For each $s \in \{2,3,4, \ldots\}$, we define the guessing number $g(G,s)$ of a graph $G$ to be the uniquely determined $\alpha$ such that $\frac{1}{s}^{n-\alpha}$ is the probability all $n$ players guess correctly their own die value (assuming the players have agreed in advance an optimal guessing strategy). The general guessing number $g(G)$ is defined as ${\rm sup}_{s=2,3,...} g(G,s)$.
\rm
\end{description}

More generally consider a pair $(f,G)$ where $f=(f_1,f_2, \ldots,f_n)$ is a fixed function with $f:A^n \rightarrow A^n$ and $G=(V,E)$ is a directed graph (not necessarily acyclic, and possibly with self loops) with $V=\{1,2,\ldots,n\}$. The values $x_1,x_2, \ldots,x_n \in A$ are selected randomly and the task for player $j$ is to guess the value of $f_j(x_1,x_2, \ldots,x_n)$. The task for the players is to find a general (deterministic) strategy that maximizes the probability that all players simultaneously guess correctly.  The {\it guessing number $g(G,f,A)$ of $(f,G)$ is the uniquely determined $\alpha$ such that $\frac{1}{s}^{n-\alpha}$ is the probability each player $j$ guess correctly the value
$f_j(x_1,x_2, \ldots,x_n)$ (assuming the players have agreed in advance an optimal guessing strategy).}
In this paper we consider in general  only the case where $f$ is the identity map, and we 
write $g(G,s)$ instead of $g(G,id,s)$.

\section{The Entropy of a graph}
\rm 
One way of thinking about the entropy of a graph informally, is to view the graph $G$ as representing a composite physical system $\cal{F}_{\rm G}$ where the nodes represent identical physical subsystems that can each be in some state $s \in A$. The (directed) edges in the graph indicate possible causal influences. The (dynamic) behavior of the composed physical system is determined by specific functions assigned to the vertices representing the underlying physical laws. Like in physical systems subject to the laws of thermodynamics, the behavior of the functions is such that the overall entropy $H_0$ of the system $\cal{F}_{\rm G}$ is maximized. Informally, the Entropy of $G$ is defined as the maximal entropy $H_0$ possible for the system $\cal{F}_{\rm G}$ when the system is subject to the causal constraints indicated by $G$ (the system is, as is assumed in thermodynamics, always in equilibrium i.e. all constraints are satisfied at any given moment).  

More formally, let $G=(V,E)$ be a (directed) graph with the set $V$ of nodes ($n:=|V|$) and set $E \subseteq V \times V$ of edges. Assume that each node $j\in V$ in $G$ has assigned a stochastic variable $x_j$ selected from some finite alphabet (or state space) $A$ with $s \geq 2$ elements. 
For each probability distribution $p$ on tuples $(x_1,x_2,\ldots,x_n) \in A^n$ we define an entropy function $H_{p}$ such that for each subset $S \subseteq \{1,2,\ldots, n\}$  the real number $H_{p}(S)$ is given by:

\begin{equation}
H_{p}(S):=\sum_{v \in A^n} p(S,v) \log_{s}(\frac{1}{p(S,v)})
\end{equation}
where $p(S,v)$ for $v=(v_1,v_2, \ldots,v_n) \in A^n$ is the probability that a tuple $(x_1,x_2,\ldots,x_n)  \in A^n$  is selected with 
$x_{s_1}=v_{s_1}, x_{s_2}=v_{s_2}, \ldots, x_{s_u}=v_{s_u}$ where $S=\{s_1,s_2,\ldots,s_u\}$.
Let $H$ denote any such entropy function $H_p$. 

For each vertex $j$ in $G$ we introduce the information equation: 

\begin{equation}
H(j | i_1,i_2,\ldots, i_d)=0 
\end{equation}
where $(i_1,j),(i_2,j),\ldots, (i_d,j) \in E$ is the edges with head 
$j$.  We refer to these $n$ information equations as the {\it information constraints} determined by $G$.
Notice that the use of logarithm in base $s$ in the definition ensures that the Entropy function $H$ is normalized. Thus it satisfies $H(j) \leq 1$ for each $j=1,2,\ldots,n$. 

Informally, the equation states that there is no uncertainty of the value of the variable $x_j$ corresponding to vertex $j$ if we are given the values of all the stochastic variables associated with  the
predecessor vertices of $j$. 

\begin{description}

\item[Definition]:

\noindent
{\it The (private) entropy $E(G,s)$ of a graph $G$ over an alphabet $A$ of size $s \in \{2,3,4, \ldots\}$ is the supremum of $H_p(1,2, \ldots, n)$ of all entropy functions $H_p$ over $A$  that satisfies the $n$ information constraints determined by $G$. 
The general entropy $E(G)$ (or just entropy) of a graph $G$ is the supremum of the entropy $E(G,s)$
for $s=2,3,4,\ldots$.}
\end{description}

As usual the conditional entropy is defined as 
$H(X | Y):=H(X,Y)-H(Y)$ and the mutual information between $X$ and $Y$ is defined as
$I(X;Y):=H(X)+H(Y)-H(X,Y)$ 
where $H(X,Y)$ is shorthand for $H(X \cup Y)$. In general we drop set clauses when possible and for example write $H(i_1,i_2,...,i_r)$ instead of $H(\{i_1,i_2,...,i_r\})$. We have $H(\emptyset)=0$.   

An entropy  function $H$ has $H(X )=0$ if and only if  $X$ is uniquely determined, and has $H(X | Y)=0$ if and only $X$ is a function $Y$ (see \cite{YoungBook} for such and other basic properties of entropy functions). Shannon showed that entropy functions satisfy a number of information inequalities 
e.g. $0 \leq H(X,Y)$ $\quad H(X,Y) \leq H(X)+H(Y)$, $\quad H(X,Y | Z) \leq H(X | Z)+H(Y | Z)$ and
$H(X | Y,Z) \leq H(X | Y) \leq H(X,Z | Y)$. 
These inequalities can all be obtained as a special case of Shannon's famous information inequality 

\begin{equation}
H(X,Y,Z) + H(Z) \leq H(X,Z) + H(Y,Z) 
\end{equation}
We will show:

\begin{description}

\item [Theorem 1]:
\it
\noindent
For each directed  graph $G$ and for each $s \in \{2,3,4, \ldots, \}$ the guessing number equals the entropy (i.e. $E(G,s)=g(G,s)$).

Furthermore, the general entropy of a graph $G$ is identical to the general guessing number of $G$
(i.e. $E(G)=g(G)$).
\rm
\end{description}

This result allows us to calculate the Entropy of a graph using quite different methods.
One can calculate lower bounds on the entropy, by constructing an explicit entropy function, but  alternatively lower bounds can be obtained by providing a guessing strategy. For all (small) graphs we tested, and we believe most graphs in general, good upper bounds on the entropy can be calculated using Shannons Information inequalities or computationally more efficiently by use of the polymatoidal axioms we will introduce in the next section. 
In general - for the vast majority of (small) graphs we analyzed by computer- it is possible to find  matching lower and upper bounds on $E(G,s)$ in this fashion.

The main reason we do not merge the Guessing number and the (private) Entropy into one concept is that we, as already mentioned,wants to introduce a number of concepts related to the Entropy of a graph  (e.g.the Shannon-entropy, the ZY-Entropy and the DFZ-entropy of a graph) and already in 
\cite{Riis06, Riis07} have defined a number of distinct types of guessing numbers of a graph (e.g. the linear guessing number and the scalar linear guessing number). 
Joint with Danchev we also considered other types of guessing numbers defined using analogues to various well known routing protocols (e.g. fractional routing \cite{Cdfz}).

\section{The S/ZY/DFZ graph entropy}

The entropy function $H$ has a number of basic properties that are very useful in computing the entropy of a graph.  In general there exist functions that satisfy all of Shannon information inequalities, but are not genuine entropy functions. An entropy-like function $f$ on an $n$-element set $V=\{1,2,\ldots,n\}$ is a map $f$ from the subsets of $V$ into $R$. \footnote{as usual we use $R$ to denote the real numbers}  
We require that an entropy-like function has $f(\emptyset)=0$. For subsets $X,Y \subseteq V$ we let $f(X,Y):=f(X \cup Y)$ and we define $f(X | Y)$ as $f(X | Y):=f(X,Y)-f(Y)$. 

We say an entropy-like function $f$ is a {\it shannon} entropy function if for $X,Y,Z \subseteq V$
we have $f(Z) + f(X,Y,Z) \leq f(X,Z) + f(Y,Z)$

As pointed out in \cite{DFZ07} this information inequality (combined with $H(\emptyset)=0$),
is known to be equivalent to the so-called polymatroidal axioms (for a map $f: 2^V \rightarrow R$): 

\medskip

\noindent
(i) $\quad$  $f(\emptyset)=0$ 

\smallskip

\noindent
(ii) $\quad$ for $X,Y \subseteq V, \quad f(X) \leq f(Y)$. 

\smallskip

\noindent
(iii) $\quad$ $f(X \cap Y) + f(X \cup Y) \leq f(X) +f(Y)$. 

\medskip

For an entropy-like function $f: P(V) \rightarrow R$ we define the mutual information $I_f$ such that

\begin{equation}
I_f(A;B | C):=f(A,C)+f(B,C)-f(A,B,C)-f(C)
\end{equation}
The special case where $C=\emptyset$,  gives $ l_f(A;B):=f(A)+f(B)-f(A,B)$.

We say an entropy-like function is a {\it Zhang-Young} entropy function (or just ZY-entropy function) if for $A,B,C,D \subseteq V$ we have 
\begin{equation}
2I_f(C;D)\leq I_f(A;B)+I_f(A; C\cup D)+3I_f(C; D | A) + I_f(C;D | B)
\end{equation}
It can be shown that each ZY-entropy function is a Shannon entropy function. Expressed in our terminology, Zhang and Young discovered \cite{ZY97} that while each entropy function is a ZY-entropy function, there exist Shannon entropy functions that fail to be ZY-entropy functions.

The condition required for functions to be a ZY-entropy function can be written as
\begin{equation}
U \leq W
\end{equation}
where

\begin{equation}
U:=2f(C)+2f(D)+f(A)+f(A \cup B)+4f(A \cup C \cup D) + f(B \cup C \cup D)
\end{equation}
and

\begin{equation}
W:=3f(C \cup D)+3f(A \cup C)+3f(A \cup D)+f(B \cup C)+f(B \cup D)
\end{equation}
Equality $U=W$ happens in a number of cases e.g. when $A=B=C=D$, when $A,B,C$ and $D$ are independent sets with respect to 
$f$ (i.e. $f(X \cup Y)=f(X)+f(Y)$  as well as when $X,Y \in \{A,B,C,D\}$ are distinct.

In \cite{DFZ07} a number of new shannon information inequalities were presented and showed to be independent. We refer to such information inequalities (that can not be derived from the ZY-information inequalities) as DFZ-information inequalities. 

Somewhat informally (depending on which of these non-shannon information inequalities we consider) we say that a function $f$ is a DFZ-entropy-like function if $I_f$ satisfies the ZY-information inequality as well the DFZ-information inequalities. 

\begin{description}

\item[Definition]:

\noindent
{\it Let $G$ be a graph. We define the  {\it S-entropy}/{\it ZY-entropy}/{\it DFZ-entropy} of $G$ as the maximal value of 
$f(\{1,2,\ldots, n\})$ for any  S/ZY/DFZ entropy-like function $f$ that satisfies the entropy constraints determined by $G$.

For a graph $G$ we let $E_{S}(G)$, $E_{ZY}(G)$ and $E_{DFZ}(G)$ denote the S-entropy resp. ZY-entrop and  DFZ-entropy of $G$. }

\end{description}

\begin{description}
\item[Proposition 2]:
\it
\noindent
For each graph $G$ and each $s \in \{2,3,4, \ldots, \}$,
$g(G,s)=E(G,s) \leq E_{DFZ} \leq E_{ZY}(G) \leq E_{S}(G)$.
Furthermore the S-Entropy $E_S(G)$ is bounded from above with the number $|V|$ 
of vertices minus the acyclic independence number of $G$. 
\end{description} 

\noindent
{\bf Proof:} The equality $g(G,s)=E(G,s)$ follows from Theorem 1. The inequality $E(G,s) 
\leq E_{DFZ}(G)$ follows from the fact that each entropy function is an DFZ-entropy-like function.
$E_{DFZ}(G) \leq E_{ZY}(G)$ follows from the fact that part of the requirement of being a DFZ-entropy-like function is that it is ZY-entropy-like.  The ZY-information inequalities
are known to imply Shannon's information inequalities, which ensures that $E_{ZY}(G) \leq E_S(G)$.
For each S-entropy-like function $f$ that satisfies the information constraints for $G$, for any subset $B$ of vertices  for each vertex $j$ with all tails in $B$,  $f(B)=f(j, B)$. Let $A \subseteq V$ is any subset with the induced graph on $V \setminus C$ being acyclic. Without loss of generality, we can assume that
the elements in $V \setminus C$ are the elements $\{1,2,\ldots,r\}$ and that there is no edge from
any vertex $i$ to a vertex $j$ if $i,j \in \{1,2, \ldots,r\}$ and $i>j$. But, then $f(C)=f(C,1)=f(C,1,2)= 
\ldots, f(V)$ and thus $f(V)=f(C) \leq |C|$. The number $|C|$ is at most $|V|$ minus the acyclic
independence number of $G$.
\hfill $\clubsuit$

\section{The Entropy of the pentagon}

It is not hard to show that $K_n$ has entropy as well as S-entropy $n-1$. The entropy (S-entropy) of an acyclic graph is $0$ (recall that the entropy $H(j)$ of a node with in-degree $0$ vanish since $H(j)=H(j | \emptyset) =0$). In general the graph $C_n$ oriented as a loop has entropy (S-entropy) $1$.
To illustrate the concepts with a less trivial example consider the pentagon $C_5$ where all edges are 
bi-directed.

The Shannon Capacity - a famous graph parameter that is notoriously difficult to calculate - of $C_5$  was eventually shown to be $\sqrt{5}=2.23..$.    We will show that pentagon $C_5$ have Entropy (as well as S-entropy, ZY-entropy  and DFZ-entropy)  $E(C_5)=2.5$ \footnote{Joint with Dantchev we showed that $C_{k}$ in general have (Shannon) Entropy  $\frac{k}{2}$ for any $k \geq 4$ }.  It turns out that the concrete Entropy over an alphabet $A$ is $2.5$ if and only if $s=|A| \in \{2,3,4,\ldots\}$ is a square number. The entropy of $C_5$ over an alphabet of two elements is  $\log_2(5)=2.32...$.

 First we show that
$C_5$ has S-entropy $\leq 2.5$ i.e that any S-entropy-like function $f$ that satisfies the information constraints determined by $C_5$ has $f(1,2,3,4,5) \leq 2.5$.

Assume that $f$ satisfies Shannon Information inequality 
 $f(A,B,C)+f(C) \leq  f(A,C) + f(B,C)$.  If we let $A=\{1\}, B=\{3\}, C=\{4,5\}$ we get
$f(1,2,3,4,5)+f(4,5)=f(1,3,4,5)+f(4,5)  \leq
 f(1,4,5)+f(3,4,5)=f(1,4)+f(3,5)  \leq f(1)+f(3)+f(4)+f(5)$.
Thus $f(1,2,3,4,5) \leq  f(1)+f(3)+f(4)+f(5)-f(4,5)$.

Next notice that  $f(1,2,3,4,5)-f(2,5) \leq f(3,4 | 2,5)=f(4| 2,5) \leq f(4|5)=
 f(4,5)-f(5)$.

Thus $f(1,2,3,4,5) \leq f(2,5)+f(4,5)-f(5) \leq f(4,5)+f(2)$.
Adding the two inequalities we get: $2f(1,2,3,4,5) \leq f(1)+f(2)+f(3)+f(4)+f(5) \leq 5$ which show that the S-entropy (and thus the entropy, as well as the ZY-entropy and DFZ-entropy ) of $C_5$ is at most $2.5$.

Next we show that the entropy of $C_n$ is at least $2.5$ if $s$ the size $s$ of the underlying alphabet 
is a square number i.e $s=t^2$ for some $t \in \{2,3,4,\ldots\}$. According to Theorem 1 (the easy direction that $g(G,s) \leq E(G,s)$) it suffices to show that $C_5$ has guessing number $\geq 2.5$. 
Assume that in each of the five nodes of $C_5$ is a player with two dice each with $t$-sides. Each of the players has access to the dice values of their two neighbors. Assume that each player's dice are labeled $L$  and $R$ and they use a guessing strategy where each player assumes that their $L$-die has the same value as the left-hand neighbors $R$-die,  and that their $R$-die has the same value as the right-hand neighbors $L$-die. Notice that this protocol divides the $10$ dice into $5$ pairs with one $L$ and one $R$ die in each pair. Notice that the players are right exactly if each of the $5$ pairs of matched $L$ and $R$-pairs have identical dice-values. This happens with probability $\frac{1}{t}^5=\frac{1}{s}^{2.5}=\frac{1}{s}^{5-2.5}$. This shows that the Entropy (and thus the S-entropy, the ZY-entropy and DFZ-entropy) of the pentagon (with bi-directed edges) is at least $2.5$. 

To show that the entropy can only be achieved when $s=|A|$ is a square number, we use the non-trivial direction of Theorem 1 ($E(G,s) \leq g(G,s)$) and show that there is no guessing strategy that achieve guessing number $2.5$ unless $s$ is a square number.  Assume the players fix a (deterministic) guessing strategy. There are $S^5$ different assignments of values $(x_1,x_2, \ldots,x_n) \in A^n$ so if we let $m$ denote the number of assignments where all players guess correctly, the guessing number of $2.5$ is achieved exactly when $m=s^{2.5}$. Since $m$ is an integer
the identity is only possible when $s$ is a square number.

Finally assume that the alphabet has size $2$. In general the guessing number of a graph $G=(V,E)$ over an alphabet of size $s$ is of the form $\log_{s}(m)$ for some $m\in \{0,1,2, \ldots, s^{|V|}\}$. Thus,  since the guessing number of $\log_2(6) >2.5$ is impossible to obtain, it suffices to show that there is a guessing strategy where all players can correctly guess their own assigned coin (coin rather than die since $s=2$) value with probability $\frac{5}{32}$. One such strategy is that each player assumes that the assigned coin values do not have three consecutive 0's and not two consecutive 1's. This condition  is satisfied globally exactly when one of the $5$ configurations occurs: 
$00101, 10010, 01001,10100, 01010$. Notice that whenever this condition is satisfied "globally", each player can (based on "local" information) deduce the value of their own coin value.
 \hfill $\clubsuit$
 
 \section{Proof of Theorem 1}
A guessing strategy $\cal{G}$  induces a map $\Psi_{\cal{G}}$ that maps $A^n$ to $A^n$. 
More specifically, $\Psi_{\cal{G}}$ maps each dice assignment
to the guess made by the individual nodes, and thus maps each tuple
$(x_1,x_2,\ldots,x_n) \in A^n$
to a tuple \newline $(z_1,z_2,\ldots,z_n) \in A^n$ representing the players' guess.

All players guess correctly their own die value exactly when the assignment $(x_1,x_2,\ldots,x_n)$
of dice values is a fix-point for $\Psi_{\cal{G}}$. From this we notice:

\noindent
{\bf Observation:} 
{\it For each guessing strategy $\cal{G}$ let $\cal{C}=\cal{C}(\cal{G}) \subseteq {\rm A}^{\rm n}$ denote the set of fixpoints for the map $\Psi_{\cal{G}} :A^n \rightarrow A^n$. The guessing strategy $\cal{G}$ corresponds to guessing number $\log_s(m)$, where $m$ is the number of fixpoints 
(i.e. $m=|\cal{C} |$).}

The set  $\cal{C} \subset {\rm A}^{\rm n}$ of fixpoints of $\Psi_{\cal{G}}$ form a code with $m$ codewords from $A^n$.  Let $p=p(\cal{G})$ be the probability distribution on $A^n$ where 
each code word in $\cal{C}$ has the same probability (in fact probability $\frac{1}{m}$) and each word in $A^n \setminus \cal{C}$ has probability $0$. Since each probability distribution $p$ on $A^n$ defines an entropy function $H_p$, this especially applies to the probability distribution $p(\cal{G})$, that therefore defines an entropy function $H_{p(\cal{G})}$.

\noindent
{\bf claim:}
{\it If $\cal{G}$ is a guessing strategy for the graph $G$, the corresponding entropy function 
$H_{p(\cal{G})}$ satisfies the information constraints determined by $G$.}

To see this consider a given vertex $v$ with predecessor vertices $w_1,w_2,...,w_d$. The guessing strategy $\cal{G}$ has the guess $z_v$  determined as a function of \newline
$x_{w_{1}},x_{w_{2}}, \ldots, x_{w_{d}}$.  
Thus each word $(c_1,c_2,.... ,c_n) \in \cal{C}$ that corresponds to a word where each player guess correctly their own assigned value, have $c_{v}$ being a function of 
$c_{w_{1}},c_{w_{2}}, \ldots, c_{w_{d}}$.  Thus the entropy of $x_{v}$ given  
$x_{w_1},x_{w_2}, \ldots, x_{w_d}$ is $0$.
In other words  \newline $H(v | w_1, w_2, \ldots, w_d)=0$.
This shows that $g(G,s) \leq E(G,s)$. 

The entropy functions that arise from guessing strategies
are constant on their support \footnote{the support of $p$ is defined to be the set of elements in $A^n$ that has non-zero probability}. Most entropy functions are non-constant on their support and do not correspond to guessing strategies. Potentially such entropy functions might turn the inequality $g(G,s) \leq E(G,s)$ into a strict inequality. To prove Theorem 1 we need to show that this cannot happen i.e. 
we need to show that $E(G,s) \leq g(G,s)$.  Assume that $H_p$ is an entropy function that corresponds to a general probability distribution $p$ on $A^n$. The support of the entropy function $H_p$ that corresponds to the probability distribution $p$, forms a code $\cal{C}=\cal{C}({\rm p}) \subseteq {\rm A}^{\rm n}$. As before we notice that the code $\cal{C}$ can be used to define an entropy function $H_{\cal{C}}$ that is constant on the $\cal{C} \subseteq {\rm A}^{\rm n}$ and vanishes outside $\cal{C}$. 

If $H_p(X | Y)=0$ then $X$ is uniquely determined by $Y$.  If $q$ is a probability distribution on $A^n$
with the same support as $p$, then $H_q(X | Y)=0$. 
Or in general if a collection of information constrains each is of the form 
$H_p(X|Y)=0$ for some probability distribution $p$, then the very same constraints are satisfied for all probability distributions $q$ with the same same support as $p$. Thus since $H_p$ and $H_{\cal{C}}$
corresponds to probability distributions with the same support ($=\cal{C}$) 
the Entropy function $H_p$ satisfies the information constraints of the graph $G$ if and only if
$H_{\cal{C}}$ satisfies the information constraints of the graph $G$. 

Furthermore, a basic result in information theory states that the symmetric probability distribution on the set of support, maximizes the entropy. Thus $H_p(1,2,\ldots,n) \leq H_{\cal{C}}(1,2,\ldots,n)$.

To complete the proof it suffices to observe that if  $H_{\cal{C}}$ is an entropy function that is defined from a code $\cal{C}$ over an alphabet $A$, and satisfies the information constraints of the graph $G$, then the code $\cal{C}$ can arises from a guessing strategy $\cal{G}$  over $A$. 

This observation can be obtained as follows: The code $\cal{C}$ satisfies the constraints that the letter $c_j$ is uniquely determined by the letters $c_{i_1},c_{i_2}, \ldots, c_{i_d}$ where $i_1,i_2, \ldots,i_d$ are the nodes with head $j$.
For each player $j$ this determine a (partial) guessing function. To make to function total, choose any value for the function (it does not matter which) on each tuple $(x_{i_1}, x_{i_2}, \ldots,x_{i_d})$ that does not have 
$x_{i_1}=c_{i_1}$, $x_{i_2}=c_{i_2}, \ldots, x_{i_d}=c_{i_d}$ for some code word $(c_1,c_2, \ldots,c_n) \in \cal{C}$. If the players follow this guessing strategy, they are all correct for all assignments 
$x_1,x_2, \ldots,x_n) \in \cal{C} \subseteq A^n$ (and possibly even for more words).
This shows that $E(G,s) \leq g(G,s)$. \hfill $\clubsuit$

\section{Minimal index codes}
Assume that we play the guessing game on a graph $G=(V,E)$ and with alphabet $A$.  Assume that we have available a public information channel that can broadcast $m$ messages and for each assignment of dice values we broadcast a single message to all the players. What is smallest value of $m$ for which there exists a general protocol such that each player is always able to deduce his/her own assigned value? This concept, defined in terms of guessing games was introduced in \cite{Riis05} without knowledge of \cite{Tel98}.  As already explained we will use the terminology of \cite{Tel98} and define the shortest index code $i(G,s)$
of $G$ (wrt. to $A$) as $\log_{s}(m)$ where $s=|A|$ denotes the size of the alphabet $A$.  

Assume that a graph $G$ has a minimal index code (over an alphabet with $s$ letters) of length 
$\log_s(m)$. Then by definition it is possible to broadcast one of $m$-messages and ensure that
all players with probability $1$ work out their dice values. This is (essentially) \footnote{essentially, since if $m$ does not divide $s^n$ each of the $m$ messages are not broadcast with identical probabilities} a factor 
$s^{n-\log_{s}(m)}$ better than the probability that the players could have achieved if they disregarded the information given from $G$ and only had available a public channel broadcasting one of $m$ messages. Using the analog to the definition of the guessing number, we define the public guessing number as $n-\log_s(m)$ i.e. as $n-i(G,s)$.

It is possible to define the entropy of the length of the minimal index code of a graph $G$.  Given
a strategy $\cal{S}$ that achieve an index code with $m$ messages (i.e. it possible to broadcast one in $m$ public messages, in such a fashion that each node in $G$ to deduce their own assigned value).
The length of the code is $\log_s(m)$. We define the entropy of the code as the entropy of the set of public messages with the probability distribution by which they will be broadcasted (this is not always the symmetric probability distribution). We let $i_{entro}(G,s)$ denote the minimal entropy (calculated using the logarithm in base $s$) of any index code for $G$.
A priori  there is no reason to believe that the minimal code contains the same number of code 
words as the code with the minimal entropy (the code might contain more words but if the occur
with less uniform probabilities their entropy might be lower). 

The value of the entropy of the index code with minimal entropy can be calculated in a fashion that
can be approximated by entropy-like functions (S-entropy-like, ZY-entropy-like as well as DFZ-entropy-like). As in the definition of the graph entropy we assume that each node $j\in V$ in $G$ has assigned a stochastic variable $x_j$ selected from some finite alphabet (or state space) $A$ with $s \geq 2$ elements. 
 Let $W$ be a finite set of potential public messages. To ensure $W$ is large enough we may assume $W=A^n$. It turns out the size of $W$ - as long as it is large enough - does not affect the result. 
 For each probability vector $p$ on $A^n \times W$  
 we define an entropy function $H_p$ analogous to equation (1). More specifically,
 for each probability distribution $p$ on tuples $(x_1,x_2,\ldots,x_n,x_w) \in A^n \times W$ we define an entropy function $H_{p}$ such that for each subset $S \subseteq \{0,1,2,\ldots, n,w\}$\footnote{where 
 $w$ is a new symbol distinct from $1,2, \ldots,n$} 
the real number 
$H_{p}(S)$ is given by (1) i.e. by

\begin{equation}
H_{p}(S):=\sum_{v \in A^n \times P} p(S,v) \log_{s}(\frac{1}{p(S,v)})
\end{equation}
but where $p(S,v)$ for $v=(v_1,v_2, \ldots,v_n, v_w) \in A^n \times W$ is the probability that a tuple 
$(x_1,x_2,\ldots,x_n, x_w)  \in A^n \times W$  is selected with 
$x_{s_1}=v_{s_1}, x_{s_2}=v_{s_2}, \ldots, x_{s_u}=v_{s_u}$ if 
$S=\{s_1,s_2,\ldots,s_u\} \subseteq \{1,2,\ldots,n,w\}$. 
Let $H$ denote any such entropy function $H_p$. 

For each vertex $j$ in $G$ we introduce the information equation: 

\begin{equation}
H(j | i_1,i_2,\ldots, i_d, w)=0 
\end{equation}
where 
$(i_1,j),(i_2,j),\ldots, (i_d,j) \in E$ is the edges with head $j$. We 
also add the information equation  

\begin{equation}
H(w | 1,2,\ldots,n)=0
\end{equation}
as well as the
equation 

\begin{equation}
H(1,2,\ldots,n)=n
\end{equation}
The task is to find the entropy function $H$
that satisfies all constraints and produce the minimal value of $H(w)$ (over $A$). We refer to 
$i_{entro}(G,s)$ as the minimal entropy of an index code over an alphabet $A=\{1,2,\ldots,s\}$ associated to $G$.

\begin{description}

\item[Proposition 3]:
{\it The minimal value of $H(w)$ (over $A$ with $s=|A|$) is at most $i(G,s)$ i.e. 
$i_{entro}(G,s) \leq i(G,s)$.}
\end{description}

\noindent
{\bf Proof:} First we show that there exists an entropy function $H_p$ 
that satisfies the constraints determined by $G$ such that 
$H_p(w) \leq i(G,s)$.  Consider a strategy $\cal{S}$ where each player 
choose a function (that depend on all incoming values as well as the public 
message).  
Consider the set $\cal{C} \subseteq A^n \times W$ of all tuples 
$(x_1,x_2,...,x_n, w) \in A^n \times W$ where $w$ is the message that 
according to the strategy 
$\cal{S}$ have to be broadcast when the players have been assigned the values 
$x_1,x_2, \ldots,x_n \in A^n$. The set $\cal{C}$ contains exactly $s^n$ elemnets. 
Let $p$ be the symmetric probability 
distribution on $\cal{C}$ 
(probability of each element being $\frac{1}{s^n}$) that vanish on 
$A^n \times W \setminus \cal{C}$. Let $H_p: P(V) \rightarrow R$ be the entropy 
function corresponding to that probability distribution.  Let $m$ denote the 
number of messages $w \in W$ that can occur this way. The entropy $H_p(W)$
is at most $\log_s(m)$ since this is the highest entropy possible on a set with $m$ elements.
Since $i(G,s)=\log_s(m)$ this shows that $H(w) \leq i(G,s)$.
\hfill $\clubsuit$

Like in the definition of private entropy we consider entropy-like functions $f$.
We define {\it the S-minimal index code} $i_S(G)$  as the minimal value of 
$f(w)$ that is possible to derive using Shannon's Information Inequalities 
(or equivalently the polymatroidal axioms). We define {\it the ZY-minimal  
index code} $i_{ZY}(G)$ as the minimal value of  $f(w)$ that is possible to 
derive using the ZY-information inequalities.  The {\it DFZ-minimal index code} 
$i_{DFZ}(G)$ is defined analogous by considering the minimal value of $f(w)$ 
that is possible for a DFZ-entropy like function.

\begin{description}

\item[Proposition 4]:  
$i(G,s) \geq i_{entro}(G,s) \geq i_{DFZ}(G) \geq i_{ZY}(G) \geq i_{S}(G)$.

\end{description}

\noindent
{\bf Proof:} The inequality $i(G,s) \geq i_{intro}(G,s)$ follows from Proposition 3.
Each entropy function is an DFZ-entropy-like function which shows that
$i_{intro}(G,s) \geq i_{DFZ}(G)$. Part of the requirement of being a DFZ-entropy-like function
is that it is ZY-entropy-like and thus $i_{DFZ}(G) \geq i_{ZY}(G)$. Finally, the ZY-condition
is known to imply Shannon's information inequalities, which ensures that $i_{ZY}(G) \geq i_{S}(G)$.
\hfill $\clubsuit$

Using the analog to the definition of the private entropy, (and Proposition 3) we define the public Entropy $E^{public}(G,s)$ of $G=(V,E)$ over an alphabet with $s$ letters as $n-i_{entro}(G,s)$ and the general public entropy as ${\rm sup}_{s=2,3,\ldots}E^{public}(G,s)$. The 
public S-entropy is defined as $n-i_{S}(G)$, the public ZY-entropy as $n-i_{ZY}(G)$ and the public DFZ-entropy as $n-i_{DFZ}(G)$. With these definitions, Proposition 4 can be restated:

\begin{description}

\item[Proposition 5]: {\it $g^{public}(G,s) \leq E^{public}(G,s) \leq E^{public}_{DFZ}(G) 
\leq E^{public}_{ZY}(G) \leq E^{public}_{S}(G)$} 

\end{description}

\noindent
{\bf Example:}  The index code of $C_5$ with minimal entropy has entropy $2.5$ (i.e. $C_5$ has public entropy $5-2.5=2.5$) over alphabets of size $4,9,16,25, \ldots$.  The length $i(C_5,s)$ of the minimal index code is also $2.5$ over alphabets of size $4,9,16,25, \ldots $. 

To see this we first show that 
$i_S(C_5) \geq 2.5$ (and thus that the public $S$-entropy is $\leq 2.5$). Assume that $f$ satisfies Shannon Information inequality
 $f(A,B,C)+f(C) \leq  f(A,C) + f(B,C)$.  If we let $A=\{1\}, B=\{3\}, C=\{4,5\}$ we get
  $f(1,2,3,4,5,w)+f(4,5,w)=f(1,3,4,5,w)+$  $f(4,5,w)  \leq
 f(1,4,5,w)+$  $f(3,4,5,w)=$  $f(1,4,w)+f(3,5,w)  \leq f(1)+f(3)+f(4)+f(5)+2f(w)$. 
Thus $f(1,2,3,4,5,w) \leq  f(1)+f(3)+f(4)+f(5)-f(4,5,w)+2f(w)$.

\noindent
 Next notice that  $f(1,2,3,4,5,w)-f(2,5,w) \leq $ $f(3,4 | 2,5,w)=$ $f(4| 2,5,w)$ 
 $\leq f(4|5,w)=$
$f(4,5,w)-f(5,w)$.  Thus $f(1,2,3,4,5,w) \leq f(2,5,w)+f(4,5,w)-f(5,w) \leq $ 
$f(4,5,w)+f(2)+f(5,w)-f(5,w)=f(4,5,w)+f(2)$.
 Adding the two inequalities we get: $2f(1,2,3,4,5,w) \leq f(1)+f(2)+f(3)+f(4)+f(5)+2f(w)$.
So far the calculation was essentially identical to the calculation for the graph entropy
with one additional node assigned a variable $x_w \in W$. Now 
the special information equation for calculating the index code
allows us to deduce that $f(1,2,3,4,5,w)=f(1,2,3,4,5)=f(1)+f(2)+f(3)+f(4)+f(5)=5$ and to
conclude that $2 \times 5 \leq 5+2f(w)$ i.e. that $2.5 \leq f(w)$ and that  $i_S(C_5) \geq 2.5$.

It is not hard to adopt the guessing strategy for the guessing game on $C_5$ and for $s \in \{4,9,16,25, \ldots \}$ broadcast  the sum (modulo $\sqrt{s}$) of each of the five L-R pairs. This shows that 
$i(C_5,2) \leq 2.5$  
\hfill $\clubsuit$

\begin{description}

\item[Proposition 5A]:

$g(G,s)=E(G,s) \geq E^{public}(G,s) \geq g^{public}(G,s)$

\end{description}

\noindent
{\bf Proof:} We already showed  $g(G,s)=E(G,s)$ (Theorem 1) and 
$E^{public}(G,s) \geq g^{public}(G,s)$ (Proposition 4). It suffices to show that $g(G,s) \geq E^{public}(G,s)$. Before we show this let us observe that $g(G,s) \geq g^{public}(G,s)$ (this was also
noticed in \cite{Riis05,Riis07}). This is because 
one message $w_0$ of the $m$ possible public messages must 
occur with probability at least $\frac{1}{m}$. If the players guess using the same strategies as in the public game, but with each player pretending that the public message is $w_0$, they all win with probability $\geq \frac{1}{m}$ and $g(G,s) \geq \log_s(m)$.

We want to show that $g(G,s) \geq E^{public}(G,s)=n-i_{entro}(G,s)$.  Assume that $(p_1,p_2, p_3, \ldots,p_u)$ is a probability vector and $0 \leq q \leq 1$ and  $\Sigma_j p_j \log_s(\frac{1}{p_j}) \leq \log_s(\frac{1}{q})$. Then
there exists $j$ such that $p_j \geq q$ (just view the left hand side as a convex linear combination
of the $\log_s(p_j)$.  Thus for some $j$ we must have $\log_s(\frac{1}{p_j}) \leq \log_s(\frac{1}{q})$).
If the players assume that the public message is $w_0$, where this message happens to be broadcasted with probability at least as high as $q$ where $q$ is chosen such that $g(G,s)=n-\log_s(\frac{1}{q})$. This shows that $g(G,s) \geq E^{public}(G,s)$.
  \hfill $\clubsuit$

The private entropy $E(G,s)$ and the public guessing number \newline 
$E^{public}(G,s)$ seem to be identical for many small graphs.  This is, however, not always the case.   

\begin{description}

\item[Proposition 6]: 
{\it There exists a graph $G$ with distinct private entropy and public guessing number 
(over a fixed alphabet). More specifically \newline $g(C_5,2)=E(C_5,2)>g^{public}(C_5,2)$}
\end{description}
{\bf Proof:}  We already calculated $E(C_5,2)=\log_2(5)=2.32...$.  
Since \newline $E^{public}(C_5,2) \leq g^{public}(C_5,2)$ and since the public guessing 
number is on the form $5-\log_2(m)$ the only two serious possibilities
for the public guessing number are that 
$g^{public}(C_5,2)=5-\log_2(7)=2.19...  $ or $g^{public}(C_5,2)=2$.  
Thus $g^{public}(C_5,2) \leq 5-\log_2(7) < \log_2(5)=E(C_5,2)$ i.e. the entropy (guessing number) 
is distinct from the public guessing number. \hfill $\clubsuit$ \newline
Our computer calculations showed that many small graphs have integer 
entropy. This can (partly) be explained by combining the fact that the 
private and public entropy tend to measure very similar notions with the 
following proposition:  

\begin{description}

\item[Proposition 7]: 
{\it If the private entropy $E(G,s)$ is identical to the public guessing number 
$g^{public}(G,s)$ and $s$ is a prime, then $E(G,s)$ is an integer.} 

\end{description}
{\bf Proof:} Let $A$ be an alphabet with $s$ letters. 
The private entropy is of the form $\log_s(m_1)$ while the public guessing number 
is of the form $n-\log_s(m_2)$ for $m_1,m_2 \in \{0,1,2, \ldots\}$. If the 
private entropy is identical to the public guessing number,
we have $n=\log_s(m_1)+\log_s(m_2)=\log_s(m_1 m_2)$ i.e. $s^n=m_1 m_2$. If 
$s$ is a prime number, this is only possible if $m_1$ as well as $m_2$ is a 
power of $s$ i.e. if $\log_s(m_1)$ and $\log_s(m_2)$ are both integers.  
\hfill $\clubsuit$

\section{Graph Entropy and its link to Network Coding}
\subsection{Historical background}
Network Coding is based on the idea that  for many network topologies messages can be transmitted more efficiently using coding functions instead of being transmitted in discrete packets
(like in traditional routing). 
This basic observation is often credited in the Network Coding literature to the highly influential paper \cite{Satellite} about Distributed Source Coding for satellite communications.  It should be pointed out
that - if we disregard issues about potential applications - the idea was known to the Circuit Complexity community at least 25 years before.  Curiously in many of Valiant papers including \cite{Valorg} where Valiant anticipated the definition of the minimal index code of a graph, Valiant  considered problems that mathematically clearly presuppose the basic observation in Network Coding. 

The main result in \cite{Satellite} (see also \cite{BR1}) was to show that Network Coding for each specific unicast problem - provided the underlying alphabet is sufficiently large - never leads to "unnecessary" information bottle necks.  More specifically, they showed that as long as the min-cut capacity  between the single source node and each specific target node allows transmission (i.e. there are no obvious obstructions), simultaneous  transmission from the source node to all target nodes is always possible. 

In our view the real significance of \cite{Satellite} was to recognize and propose Network Coding as an alternative approach to data management in communication networks. Later the notion of Network Coding has been broadened and now include, ideas in wireless communication, peer-to-peer files sharing (e.g. as in Avalanche by Microsoft).  For a general introduction to Network coding see
for example \cite{intro1} or \cite{intro2}.

\subsection{Solvability and graph entropy}
A multiple unicast network problem is an acyclic graph with $k$ input nodes and $k$ output nodes. Each input is demanded exactly at one output.  Each node (except for the input nodes) is assigned a coding function and functions compose in the obvious manner. 

The acyclic independence number of a graph $G=(V,E)$ is defined as the size of the maximal induced acyclic subgraph of $G$ \footnote{the task $(G,r)$ of deciding if the input graph $G$ has independence number $\geq r$ is NP-complete \cite{bang}}  \cite{bang}.  In general if $B \subseteq V$ we can obtain a network $N$ by duplicating $B$, such that each node in $B$ is copied into two nodes, one being an input node and one being an output node.  We say $B$ is {\it a split} if the resulting network $N$
is acyclic (i.e. if the induced graph on $V \setminus B$ is acyclic). The number of nodes in $B$ is the {\it size of the split}.  {\it A minimal split} is a split of minimal size and {\it the minimal split} of a graph is the size a minimal split.  The minimal split of an acyclic graph is $0$ since the empty set can serve as a split set.  

\includegraphics[scale=0.45]{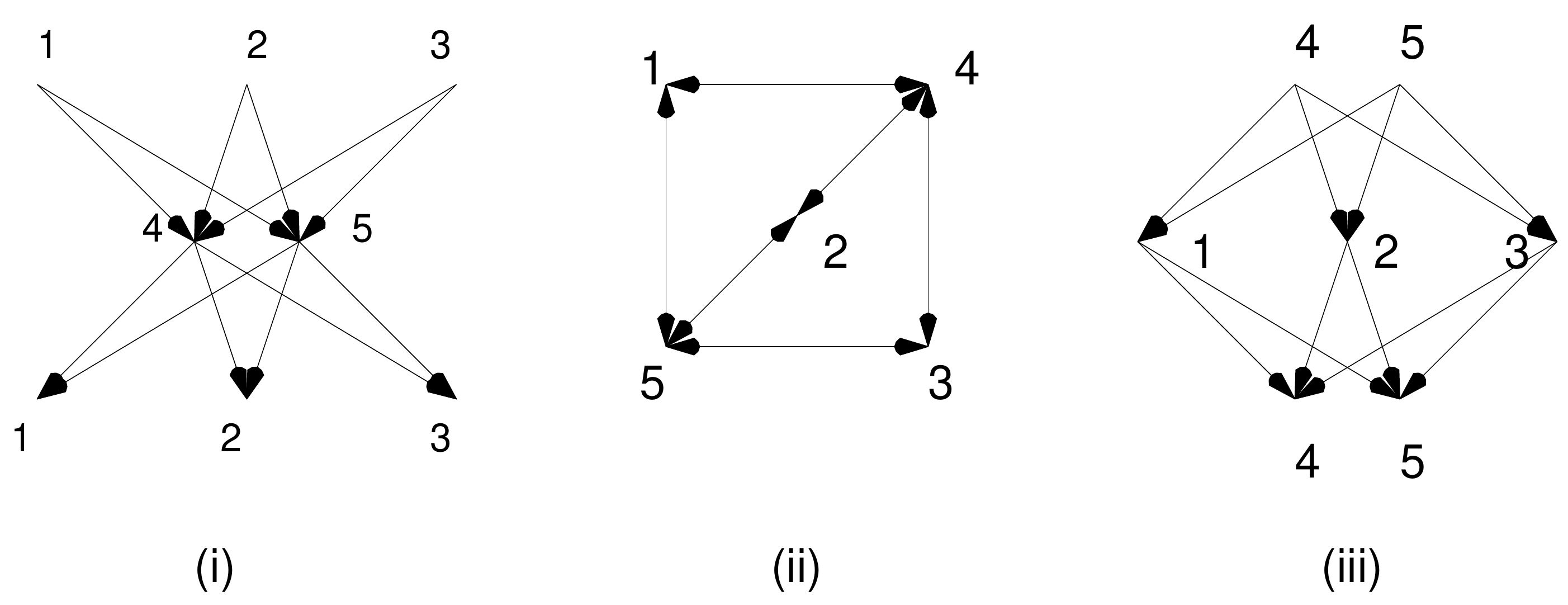} 

As an example consider the communication Network {\bf (i)} where the 
inputs in vertices 1,2 and 3 at the to
are required at vertices 1,2 and 3 at the bottom.  The network is unsolvable since all messages
have to pass through nodes 4 and 5 that create an information bottleneck. Consider the graph {\bf (ii)} that appeared by identifying each input node with its corresponding output node. If the graph {\bf (ii)}  is split with split set $\{4,5\}$ we get network {\bf (iii)}.  This network
is clearly solvable which shows that graph {\bf (ii)} has entropy $2$.

We can show (see \cite{Riis06} and \cite{Riis07} for a related result):

\begin{description}

\item[Theorem 8]:
\it
Let $G$ be a directed graph, and let $N_G$ be a communication network (circuit) that occurs by splitting
$G$ with a split of size $k$. Fix an alphabet $A$ and assume that guessing functions $g_1,g_2,\ldots, g_n$ are assigned to each node in $G$. Further, assume that the very same functions (with the arguments in a specific order that depends on the split) are assigned to the corresponding computational nodes in $N_G$. Then the following are equivalent:

\noindent
(i) The probability that  by using the guessing strategy determined by $g_1,g_2, \ldots, g_n$ all players guess correctly their own assigned values is (at least) $\frac{1}{s^{n-k}}$

\noindent
(ii) The coding functions $g_1 \circ \tau_1, g_2 \circ \tau_2, \ldots, g_n \circ \tau_n$ define a solution to the Network $N_G$ where $\tau_1,\ldots,\tau_n$ are suitable permutations of the arguments determined by the split.

\noindent
If the split is not minimal (i) and (ii) are never satisfied and the theorem becomes vacuous.

\rm
\end{description}

\noindent
{\bf Proof:} Assume that  the players use the functions $g_1,g_2, \ldots,g_n$ as guessing functions on the graph $G=(V,E)$. Let $B$ be the split. The induced graph on $V \setminus B$ is acyclic. Without loss of generality we can assume 
that the nodes $1,2, \ldots,k$  belong to $B$ while the nodes $k+1,k+2, \ldots,n$ belongs to $V \setminus B$. Assume that the very same functions $g_1,g_2, \ldots, g_n$
are assigned as coding functions in the corresponding network $N_G$. 
Let $\vec{x}=(x_1,x_2,\ldots,x_n)\in A^n$ be a random assignment of values to the nodes in $G$.
Let $(z_1(\vec{x}),z_2(\vec{x}),\ldots,z_n(\vec{x}))\in A^n$ denote the "guess" by the players determined
by the guessing functions $g_1,g_2, \ldots,g_n$. 

Condition (i) can then be stated as:
\begin{equation}
p(z_{1}=x_{1}, z_{2}=x_{2}, \ldots, z_{n}=x_{n})= (\frac{1}{s})^{n-k}
\end{equation}
and condition (ii) can be stated as:
\begin{equation}
p(z_1=x_1, z_2=x_2, \ldots, z_k=x_k | z_{k+1}=x_{k+1}, z_{k+2}=x_{k+2}, \ldots, z_{n}=x_{n})=1
\end{equation}

In general $P(U | V)= \frac{P(U \cap V)}{P(V)}$ and thus

\[p(z_1=x_1, z_2=x_2, \ldots, z_k=x_k | z_{k+1}=x_{k+1}, z_{k+2}=x_{k+2}, \ldots, z_{n}=x_{n})=\]
\begin{equation}
\frac{ p(z_{1}=x_{1}, z_{2}=x_{2}, \ldots, z_{n}=x_{n})}{ p(z_{k+1}=x_{k+1}, z_{k+2}=x_{k+2}, \ldots, z_{n}=x_{n})} 
\end{equation}

In general,  for  {\it any} choice of guessing functions (coding functions)
we have

\begin{equation}
p(z_{k+1}=x_{k+1}, z_{k+2}=x_{k+2}, \ldots, z_{n}=x_{n})= (\frac{1}{s})^{n-k}
\end{equation}
since the induced graph of $G$ restricted to the nodes $k+1,k+2, \ldots,n$ is acyclic
so the involved probabilities are independent. Substituting (16) into (15) we get that
\[ (\frac{1}{s})^{n-k} p(z_1=x_1, z_2=x_2, \ldots, z_k=x_k | z_{k+1}=x_{k+1}, z_{k+2}=x_{k+2}, \ldots, z_{n}=x_{n})\]
\begin{equation}
= p(z_{1}=x_{1}, z_{2}=x_{2}, \ldots, z_{n}=x_{n})
\end{equation}
i.e. that (13) holds if and only if (14) holds. \hfill $\clubsuit$

To illustrate the theorem, consider the Networks (b)-(d) and the assignments of guessing functions $g_1,g_2, \ldots, g_7$ for the graph {\bf a}.  Each Networks (b)-(d) appears by "splitting" the 
graph {\bf a}  in different ways. 

Consider as an example Network {\bf b} \footnote{we prefer to represent multiple unicast problems in Network coding by using the line-graph. This way  coding functions are attached to the nodes rather to the edges}. The network in {\bf b} was introduced in \cite{DFZ05, DFZ06} where the authors showed that the network
has only linear solutions over fields of characteristic $2$ (i.e when the alphabet contains $2^r$ elements for some $r \in \{1,2,3, \ldots\}$). Our concern is to illustrate Theorem 8.  Let us recall that solving network {\bf b} consists of finding coding functions $g_1,g_2,\ldots,g_7$ such that for any choice of input messages $x_1,x_5,x_2 \in A$, (for nodes 1,5 and 2 at the top) these messages can be reproduced at their corresponding output nodes (nodes 1,5 and 2 at the bottom). Since each vertex has two predecessors,  each coding function is a map $g_j: A^2 \rightarrow A$. 
The coding functions compose in the obvious way and they solve the network if for any input $(x_1,x_5,x_2)$ the output (viewing the network as a circuit) is $(x_2,x_5,x_1)$  (reading from left to right).  

\includegraphics[scale=0.50]{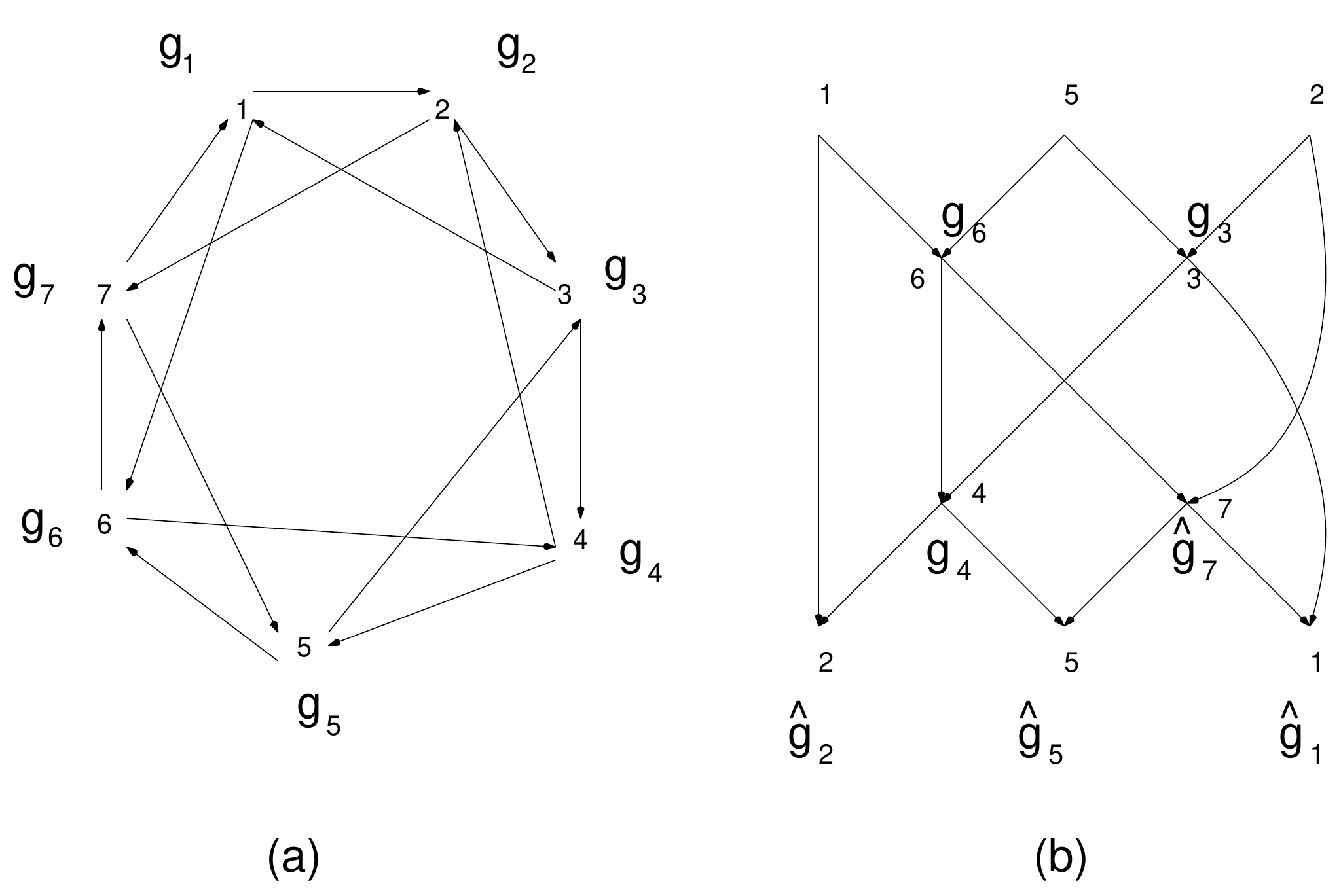} 

\includegraphics[scale=0.60]{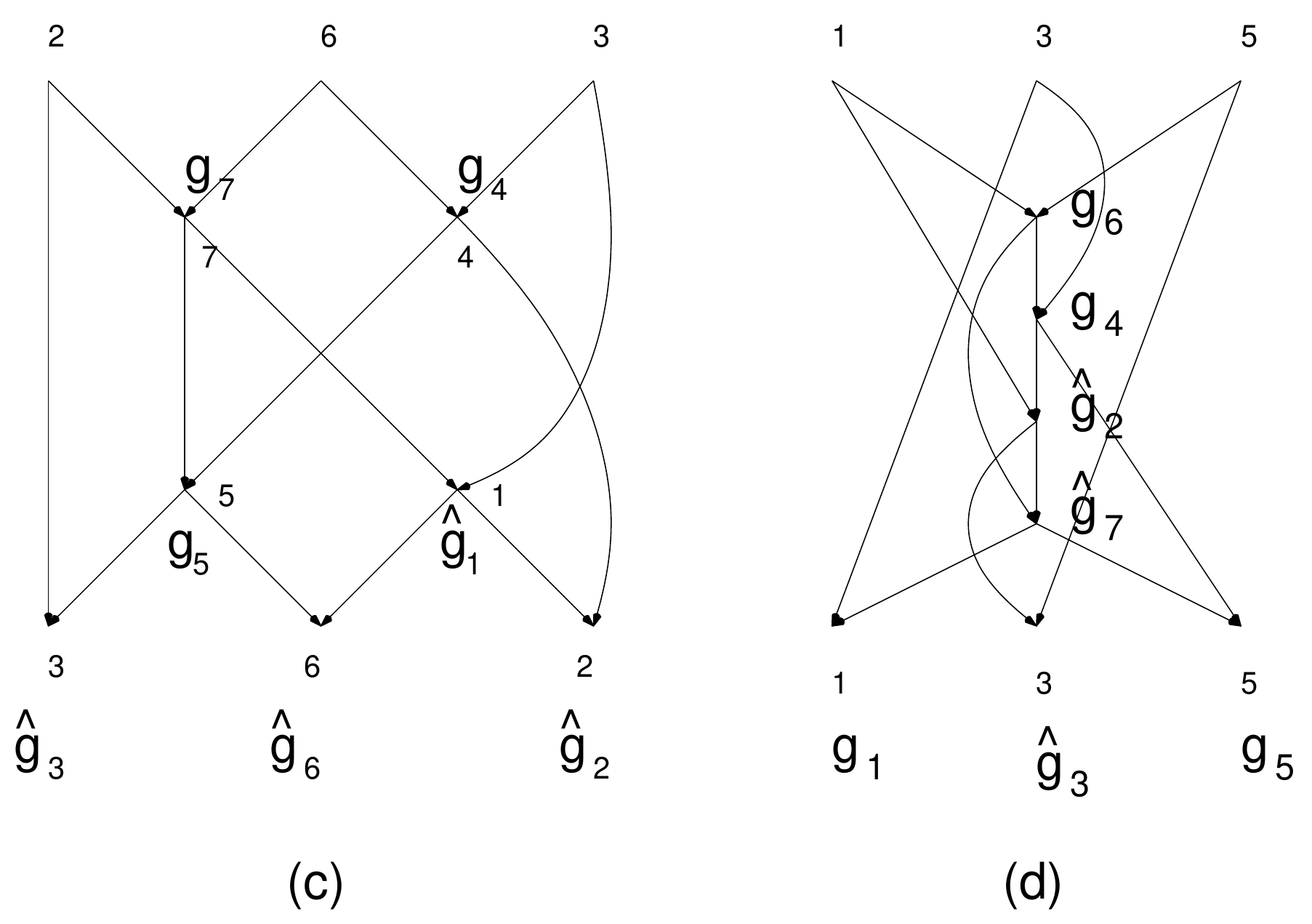} 

If we in Network {\bf b} identify the three sources at the top (nodes $1,2$ and $5$) with the corresponding output nodes at the bottom (nodes $1,2$ and $5$), we get the graph {\bf a}. The functions $g_1,g_2,\ldots,g_7$ are assigned to the corresponding nodes in Network {\bf b}.  
The spilt  $\{1,2,5\}$ leads to network {\bf b} while the split $\{2,3,6\}$ leads to network {\bf c}. Other splits e.g.   $\{3,4,7\}$ lead to a network that is isomorphic to network {\bf a} and network {\bf b} .  The split
$\{ 1,3,5\}$ leads to network {\bf d} and the split $\{2,4,6\}$ would lead to Network isomorphic to network ${\bf d}$.  We assigned guessing functions to the nodes in the graph ${\bf a}$. Notice that each split leads to different assignments of coding functions in the resulting network.
 {\it Please keep in mind that the arguments of coding functions has to be defined appropriately}.  Let 
 $\tau(u,v)=(v,u)$ and define in general  $\hat{h}:=h \circ \tau$. Then if we define the coding functions such that they take the left incoming value as first argument, and right incoming  as second argument, we get a network coding solution for Network {\bf b} by assigning the coding functions $\hat{g}_1,\hat{g}_2,g_3,g_4, \hat{g}_5, g_6$ and $\hat{g}_7$ as indicated.
%\includegraphics[scale=0.35]{N1Split3.pdf} 
%Splitting with split set $\{1,3,5\}$ leads to the network in figure(e), while the splitting with split set $\{2,4,6\}$ leads to the network in figure(f). The networks (e) and (f) are identical however the assignment of coding functions inherited from the guessing strategy differers. On the other hand it is not hard to show that network (e) and (f) not are isomorphic to networks (b),(c) and (d). 

One direct consequence of theorem 8 is that {\it the functions $g_1,g_2,\ldots,g_7$, provide 
a solution for one of the networks  $(b),(c)$ or $(d)$ if and only if they provide a solution for each of the networks ($(b),(c)$ and $(d)$).}  
Combining Theorem 1 that the entropy is identical to the guessing number, as well as Proposition 3, i.e. that the public entropy of a graph $G=(V,E)$ is identical to $|V|$ minus the minimal index code of $G$ (since the minimal index code is the very same notion as the information defect \cite{Riis05,Riis07}),  Theorem 4 in \cite{Riis07} can be rephrased as follows:

\begin{description}

\item[Corollary 9]: 
{\it A multiple unicast network $N$ with $k$ input and $k$ output nodes, is a solvable Network over an alphabet $A$ if and only if the graph $G_N$ has (private/public) entropy $k$ over $A$}
\end{description}
 
Corollary 9 is not correct if we remove reference to a specific alphabet $A$ and just consider the general Entropy. This follows essentially by considering the unsolvable network $N$ with coding capacity $1$ constructed in \cite{Cdfz}.

 \subsection{Links to Network coding capacity}
In \cite{DFZ07} the authors define the coding capacity of a communication network. Here we are only interested in the special case of coding capacity $1$ (more specifically coding capacity $(1,1)$). Let $N$ be a communication network i.e.
an acyclic graph with a set $I$ of $k$ input nodes and $k$ output nodes. We say that N has S-coding capacity $(1,1)$ (ZY-coding capacity $(1,1)$) if there exists a S-entropy-like function $f$ (ZY-entropy-like function $f$) such that

\bigskip

\noindent
{\bf (*)}  Conditions N1,N2 and N3 from \cite{DFZ07} for $(1,1)$ coding capacity:

\smallskip

\noindent
N1:   $f(I)=k$, 

\smallskip

\noindent
N2:  $f(j) \leq 1$ for each node $j$ 

\smallskip

\noindent
N3:  $f(j | i_1,i_2, \ldots,i_d)=0$ for each node $j$ where $i_1,i_2, \ldots,i_d$ are the predecessors 
of $j$. 

\bigskip

The coding capacity conditions N1,N2 and N3 (for capacity $(1,1)$  networks) differ of course in general from the definition of graph entropy. The coding capacity is defined for acyclic networks, with special input nodes, rather than for graphs. Condition N1 for calculating the coding capacity assumes that the certain nodes are assigned independent variables, while there is no such assumption in the definition of graph entropy. On the other hand, the graph entropy assumes that a certain entropy is being optimized, while there is no such assumption in the definition of the coding capacity. However, the two notions Graph entropy and Coding capacity are in some sense identical in the special case where the maximal value of $f(1,2,\ldots,n)=f(I)$ (in which case condition N1 holds).

In other words, if we disregard condition N1 and just try to maximize $f(1,2,\ldots,n)=f(I)$ like we need to calculate the graph entropy, then since the network is assumed to have S-coding capacity $(1,1)$ (ZY-coding capacity $(1,1)$) 
we conclude that $f(I) =k$ which ensures that condition N1 is automatically satisfied.  Conversely, assume that $G_N$ has S-entropy $k$ (ZY-entropy $k$). But then since 
$f(I)=f(1,2,\ldots,k)=f(1,2,\ldots, n)=k$, we conclude that the conditions N1,N2 and N3 are satisfied. Hence:

\begin{description}

\item[Proposition 10]:  
{\it A network $N$ with $k$ input-output pairs, has shannon 
coding capacity (ZY-coding capacity) 
$(1,1)$ if and only if
 the graph $G_N$ that occurs by identifying each input node with its corresponding output node
 has S-entropy (ZY-entropy) $k$ (or $\geq k$).}
\end{description}

 \begin{description}
 \item [Theorem 11]:
 \it
\noindent
There exists a graph $G$ with $E_{ZY}(G) < E_S(G)$. Thus it has a general entropy $E(G)$ that cannot be calculated using Shannon's information inequalities.  
% \noindent
%Furthermore,  $i_{ZY}(G) > i_S(G)$,  thus the optimal lower bound on the size of the minimal index code $i(G)$ cannot be calculated using only Shannon's information inequalities.   
\rm
  \end{description}
 
\noindent
{\bf Proof:}  From \cite{DFZ07} we know that the Vamos communication network $N'$ modified to the multiple-unicast-situation has S-capacity $1$, but has strictly smaller ZY-capacity (at most $12/13$). 
The network can be represented as a circuit (by passing to the line graph).  This does not change anything since the involved coding functions are the same and the resulting circuit $N$ still has
S-capacity 1 and strictly less ZY-capacity (that is at most $12/13$). We then identify each source node with its corresponding target node and obtain a graph $G_N$.  
It follows from Theorem 10, that since the S-capacity of $N$ is $(1,1)$ - the $S$-entropy of $G_N$ is $k$ (where $k$ is the number of input nodes in $N$).  
The ZY-capacity of $N$ is not $(1,1)$ - so again according to Theorem 10 the ZY-entropy of $G_N$ is strictly less than $k$   \hfill $\clubsuit$

We do not attempt to calculate the specific difference between $G$'s S-entropy and ZY-entropy, and this is left as an open problem.

\section{Final remarks and acknowledgments}
Proposition 10 and Theorem 11  as well as many concepts introduced in the  paper depends strongly on the work  by Dougherty, Freiling and Zeger (\cite{DFZ05,DFZ06,DFZ07}) as well as the
earlier work by Zhang and Young (\cite{ZY97}). The use of entropy functions for multiple-unicast networks was introduced in \cite{DFZ05}. 

\bibliographystyle{abbrv}
\bibliography{bibioEntro}  % sigproc.bib is the name of the Bibliography in this case

\begin{thebibliography}{10}

\bibitem{Alon}
N.~Alon.
\newblock The shannon capacity of a union.
\newblock {\em Combinatorica}, 18:301--310, 1998.

\bibitem{bang}
J.~Bang-Jensen and G.~Gutin.
\newblock {\em Digraphs Rheory, Algorithms and Applications}, pages 1--772.
\newblock Springer-Verlag, 2007.

\bibitem{Tel06}
Z.~Bar-Yossef, Y.~Birk, T.~Jayram, and T.~Kol.
\newblock index coding with side information.
\newblock {\em FOCS}, 2006.

\bibitem{Tel98}
Y.~Birk and T.~Kol.
\newblock Coding-on-demand by an informed source (iscod) for efficient
  broadcast of different supplemental data to caching clients.
\newblock {\em IEEE Transactions on Information Theory}, 52:2825--2830, 2006.

\bibitem{Cdfz}
J.~Cannons, R.~Dougherty, C.~Freiling, and K.~Zeger.
\newblock Network routing capacity.
\newblock {\em IEEE Trans. Inf. Theory}, 52(3), March 2006.

\bibitem{DFZ05}
R.~Dougherty, C.~Freiling, and K.~Zeger.
\newblock Insufficiency of linear coding in network information flow.
\newblock {\em IEEE Transactions on Information Theory}, 51(8):2745--2759,
  August 2005.

\bibitem{DFZ06}
R.~Dougherty, C.~Freiling, and K.~Zeger.
\newblock Unachievability of network coding capacity.
\newblock {\em IEEE Transactions on Information Theory}, 52(6):2365--2372, Jun
  2006.

\bibitem{DFZ07}
R.~Dougherty, C.~Freiling, and K.~Zeger.
\newblock Networks,matroids, and non-shannon information inequalities.
\newblock {\em IEEE Transactions on Information Theory}, 53(6):1949--1969, June
  2007.

\bibitem{DFZ06A}
R.~Dougherty and K.~Zeger.
\newblock Nonreversibility and equivalent constructions of multiple-unicast
  networks.
\newblock {\em IEEE Transactions on Information Theory}, 52(11):5067--5077, Nov
  2006.

\bibitem{intro2}
C.~Fragouli and E.~Soljanin.
\newblock {\em Network Coding Fundamentals}.
\newblock Now Publishers Inc, 2007.

\bibitem{med}
R.~Koetter, M.~Effros, T.~Ho, and M.~Medard.
\newblock On coding for non-multicast networks.
\newblock {\em 41th Annual Allerton Conference on Communication Control and
  computing}, 2003.

\bibitem{rev}
R.~Koetter, M.~Effros, T.~Ho, and M.~Medard.
\newblock Network codes as codes on graphs.
\newblock {\em In Proceeding of CISS}, 2004.

\bibitem{BR1}
R.~Koetter and M.~Medard.
\newblock An algebraic approach to network coding.
\newblock In {\em Proocedings of the 2001 IEEE International Symposium on
  Information Theory}.

\bibitem{Leh}
A.~Lehman and E.~Lehman.
\newblock Complexity classification of network information problems.
\newblock {\em 41th Annual Allerton Conference on Communication Control and
  computing}, October 2003.

\bibitem{Tel07}
E.~Lubetzky and U.~Stav.
\newblock Non-linear index coding outperform the linear optimum.
\newblock {\em STOCS}, 2007.

\bibitem{Val18}
L.~Pippenger, N.~Valiant.
\newblock Shifting graphs and their applications.
\newblock {\em JACM}, 23:423--432, 1976.

\bibitem{Pudlak97}
P.~Pudlak, V.~Rodl, and J.~Sgall.
\newblock Boolean circuits, tensor ranks, and communication complexity.
\newblock {\em SIAM Journal of Computing}, 26(3):605--633, June 1997.

\bibitem{Riis04}
S.~Riis.
\newblock Linear versus non-linear boolean functions in network flow.
\newblock In {\em 38th Annual Conference on Information Science and Systems
  (CISS), Princeton, NJ, March 2004}.

\bibitem{Riis07}
S.~Riis.
\newblock Information flows, graphs and their guessing numbers.
\newblock {\em Electronic journal of Combinatorics}, R44:1--17, 2006.

\bibitem{Riis05}
S.~Riis.
\newblock Utilising public information in network coding.
\newblock {\em General Theory of Information Transfer and Combinatorics}, pages
  861--897, 2006.

\bibitem{Riis06}
S.~Riis.
\newblock Reversible and irreversible information networks.
\newblock {\em IEEE transactions of InformationTheory}, 2007.

\bibitem{Valorg}
L.~Valiant.
\newblock Graph-theoretic arguments in low-level complexity.
\newblock In {\em Lecture Notes in Computer Science}, volume~53, pages
  162--176. Springer, 1977.

\bibitem{Val}
L.~Valiant.
\newblock Why is boolean circuit complexity theory difficult?
\newblock In M.~Pattorson, editor, {\em Springer Lecture Series}, pages 84--94,
  1992.

\bibitem{Vitnew}
J.~Vitter.
\newblock External memory algorithms and data structures: Dealing with massive
  data.
\newblock {\em ACM Computing Surveys}, 33(2):209--271, June 2001.

\bibitem{intro1}
R.~Yeung, S.~Li, N.~Cai, and Z.~Zhang.
\newblock {\em Network Coding Theory}.
\newblock Now Publishers Inc, 2006.

\bibitem{Satellite}
R.~Yeung and Z.~Zhang.
\newblock Distributed source coding for satellite communications.
\newblock {\em IEEE Trans. Inform. Theory}, (IT-45):1111--1120, 1999.

\bibitem{YoungBook}
R.~Young.
\newblock {\em A First Course in Information Theory}, pages 1--412.
\newblock Springer-Verlag, 2002.

\bibitem{ZY97}
Z.~Zhang and R.~Young.
\newblock A non-shannon type information inequality.
\newblock {\em IEEE Trans. Inf. Theory}, 43(6):1982--1985, Nov 1997.

\bibitem{ZY97A}
Z.~Zhang and R.~Young.
\newblock On characterization of entropy function via information inequalities.
\newblock {\em IEEE Trans. Inf. Theory}, 44(4):1440--1452, Jul 1998.

\end{thebibliography}

%\balancecolumns
% That's all folks!
\end{document}